\documentclass[11pt,reqno,tbtags]{article}
\usepackage{amsmath,color,amsthm}

\numberwithin{equation}{section}

\newtheorem{theorem}{Theorem}[section]
\newtheorem{corollary}[theorem]{Corollary}
\newtheorem{lemma}[theorem]{Lemma}
\newtheorem{proposition}[theorem]{Proposition}
\newtheorem{claim}[theorem]{Claim}
\newtheorem{example}[theorem]{\sl Example}

\theoremstyle{definition}
\newtheorem{remark}[theorem]{Remark}

\newcommand{\lf}{\left\lfloor}

\newcommand{\rf}{\right\rfloor}

\newcommand{\EE}{{\bf  E}}

\newcommand{\PP}{{\bf  P}}

\newcommand{\diag}{{\rm diag}}

\newcommand{\xh}{\hat{x}}

\newcommand{\ph}{\hat{p}}
\newcommand{\Ph}{\widehat{P}}
\newcommand{\Gh}{\widehat{G}}
\newcommand{\Xh}{\widehat{X}}
\newcommand{\Th}{\widehat{T}}
\newcommand{\Vh}{\widehat{V}}
\newcommand{\gLh}{\widehat{\gL}}

\newcommand{\begp}{\begin{proposition}}
\newcommand{\enp}{\end{proposition}}
\newcommand{\begt}{\begin{theorem}}
\newcommand{\ent}{\end{theorem}}
\newcommand{\begl}{\begin{lemma}}
\newcommand{\enl}{\end{lemma}}
\newcommand{\begc}{\begin{corollary}}
\newcommand{\enc}{\end{corollary}}
\newcommand{\begcl}{\begin{claim}}
\newcommand{\encl}{\end{claim}}
\newcommand{\begr}{\begin{remark}}
\newcommand{\enr}{\end{remark}}
\newcommand{\begal}{\begin{algorithm}}
\newcommand{\enal}{\end{algorithm}}
\newcommand{\begd}{\begin{definition}}
\newcommand{\enf}{\end{definition}}
\newcommand{\begx}{\begin{example}}
\newcommand{\enx}{\end{example}}
\newcommand{\bega}{\begin{array}}
\newcommand{\ena}{\end{array}}

\newcommand{\ignore}[1]{}

\def\rompar(#1){\textup(#1\textup)}    

\newcommand\eps{\varepsilon}

\newcommand\gl{\lambda}
\newcommand\gL{\Lambda}

\newcommand{\refS}[1]{Section~\ref{#1}}
\newcommand{\refT}[1]{Theorem~\ref{#1}}

\newcommand{\refL}[1]{Lemma~\ref{#1}}
\newcommand{\refR}[1]{Remark~\ref{#1}}

\newcommand\noqed{\renewcommand{\qed}{}} 

\newcommand\Roesler{R\"{o}sler}

\begin{document}

\setcounter{page}{0}
\thispagestyle{empty}

\begin{center}
{\Large \bf The passage time distribution for a birth-and-death chain:\ Strong stationary duality gives a first stochastic proof \\}
\normalsize

\vspace{4ex}
{\sc James Allen Fill\footnotemark} \\
\vspace{.1in}
Department of Applied Mathematics and Statistics \\
\vspace{.1in}
The Johns Hopkins University \\
\vspace{.1in}
{\tt jimfill@jhu.edu} and {\tt http://www.ams.jhu.edu/\~{}fill/} \\
\end{center}
\vspace{3ex}

\begin{center}
{\sl ABSTRACT} \\
\end{center}

A well-known theorem usually attributed to Keilson states that, for an irreducible continuous-time birth-and-death chain on the nonnegative integers and any~$d$, the passage time from state~$0$ to state~$d$ is distributed as a sum of~$d$ independent exponential random variables.  Until now, no probabilistic proof of the theorem has been known.  In this paper we use the theory of strong stationary duality to give a stochastic proof of a similar result for discrete-time birth-and-death chains and geometric random variables, and the continuous-time result (which can also be given a direct stochastic proof) then follows immediately.  In both cases we link the parameters of the distributions to eigenvalue information about the chain.
We also discuss how the continuous-time result leads to a proof of the Ray--Knight theorem.    

Intimately related to the passage-time theorem is a theorem of Fill that any fastest strong stationary time~$T$ for an ergodic birth-and-death chain on $\{0, \dots, d\}$ in continuous time with generator~$G$, started in state~$0$, is distributed as a sum of~$d$ independent exponential random variables whose rate parameters are the nonzero eigenvalues of $- G$.  Our approach yields the first (sample-path) construction of such a~$T$ for which individual such exponentials summing to~$T$ can be explicitly identified.

\bigskip
\bigskip

\begin{small}

\par\noindent
{\em AMS\/} 2000 {\em subject classifications.\/}  Primary 60J25;
secondary 60J35, 60J10, 60G40.
\medskip
\par\noindent
\emph{Date.} Revised May~18, 2009. 
\end{small}

\footnotetext[1]{Research supported by NSF grant DMS--0406104,
and by The Johns Hopkins University's Acheson~J.\ Duncan Fund for the
Advancement of Research in Statistics.}

\newpage

\section{Introduction and summary}
\label{S:intro}

A well-known theorem usually attributed to Keilson~\cite{Keil} (Theorem~5.1A, together with Remark 5.1B; see also Section~1 of~\cite{Keillog}), but which---as pointed out by Laurent Saloff-Coste via Diaconis and Miclo~\cite{DM}---can be traced back at least as far as Karlin and McGregor~\cite[equation~(45)]{KM}, states that, for an irreducible continuous-time birth-and-death chain on the nonnegative integers and any~$d$, the passage time from state~$0$ to state~$d$ is distributed as a sum of $d$ independent exponential random variables with distinct rate parameters.  Keilson, like Karlin and McGregor, proves this result by analytical (non-probabilistic) means.

Modulo the distinctness of the rates, and with additional information (see, e.g., \cite{BS}) relating the exponential rates to spectral information about the chain, the theorem can be recast as follows. 

\begin{theorem}\label{T:Keil}
Consider a continuous-time birth-and-death chain with generator~$G^*$ on the state space $\{0, \dots, d\}$ started at~$0$, suppose that~$d$ is an absorbing state, and suppose that the other birth rates $\gl^*_i$, $0 \leq i \leq d - 1$, and death rates $\mu^*_i$, $1 \leq i \leq d - 1$, are positive.  Then the absorption time in state~$d$ is distributed as the sum of~$d$ independent exponential random variables whose rate parameters are the~$d$ nonzero eigenvalues of $- G^*$.
\end{theorem}
 
There is an analogue for discrete time:
\begin{theorem}\label{T:Keildisc}
Consider a discrete-time birth-and-death chain with transition kernel~$P^*$ on the state space $\{0, \dots, d\}$ started at~$0$, suppose that~$d$ is an absorbing state, and suppose that the other birth probabilities $p^*_i$, $0 \leq i \leq d - 1$, and death probabilities $q^*_i$, $1 \leq i \leq d - 1$, are positive.  Then the absorption time in state~$d$ has probability generating function
$$
u \mapsto \prod_{j = 0}^{d - 1} \left[ \frac{(1 - \theta_j) u}{1 - \theta_j u} \right],
$$
where $- 1 \leq \theta_j < 1$ are the~$d$ non-unit eigenvalues of $P^*$.
\end{theorem}

In this paper we will give a stochastic proof of \refT{T:Keildisc} under the additional hypothesis that all eigenvalues of~$P^*$ are (strictly) positive; as we shall see later (\refL{L:mono}), this implies another condition key to our development, namely, that
\begin{equation}
\label{mono+}
p^*_{i - 1} + q^*_i < 1, \qquad 1 \leq i \leq d.
\end{equation}
Whenever~$P^*$ has nonnegative eigenvalues, the conclusion of \refT{T:Keildisc} simplifies:
\medskip

\noindent
\emph{The absorption time in state~$d$ is distributed as the sum of $d$ independent geometric random variables whose failure probabilities are the non-unit eigenvalues of $P^*$.}
\medskip

The special-case of \refT{T:Keildisc} for positive eigenvalues establishes the theorem in general by the following argument (which is unusual, in that it is not often easy to relate characteristics of a chain to a ``lazy'' modification).  Choose any $\eps \in (0, 1/2)$ and apply the special case of \refT{T:Keildisc} to the ``lazy'' kernel $P^*(\eps) := (1 - \eps) I + \eps P^*$.  Let~$T^*$ and $T^*(\eps)$ denote the respective absorption times for~$P^*$ and $P^*(\eps)$.  Then $T^*(\eps)$ has probability generating function (pgf)
\begin{equation}
\label{pgf1}
\EE\,s^{T^*(\eps)} = \prod_{j = 0}^{d - 1} \left[ \frac{\eps (1 - \theta_j) s}{1 - (1 - \eps (1 - \theta_j)) s} \right].
\end{equation}
But the conditional distribution of $T^*(\eps)$ given~$T^*$ is negative binomial with parameters~$T^*$ and~$\eps$, so the pgf of $T^*(\eps)$ can also be computed as
\begin{equation}
\label{pgf2}
\EE\,s^{T^*(\eps)} = \EE\,\EE\left(\left.s^{T^*(\eps)}\right|T^*\right) = \EE\left(\frac{\eps s}{1 - (1 - \eps) s}\right)^{T^*}.
\end{equation}
Equating~\eqref{pgf1} and~\eqref{pgf2} and letting $u := \eps s / [1 - (1 - \eps) s]$, we find, as desired,
$$
\EE\,u^{T^*} = \left[ \frac{(1 - \theta_j) u}{1 - \theta_j u} \right].
$$
\smallskip

Later in this section we explain how \refT{T:Keil} follows from \refT{T:Keildisc}, but in \refS{S:cont} we will also outline a direct stochastic proof of \refT{T:Keil}.

\begin{remark}
\label{discconds}
(a) Theorems~\ref{T:Keil} and~\ref{T:Keildisc} are the starting point of an in-depth consideration of separation cut-off for birth-and-death chains in~\cite{DS}.

(b) By a simple perturbation argument, Theorems~\ref{T:Keil} and~\ref{T:Keildisc} extend to all birth-and-death chains for which the birth rates $\gl^*_i$ (respectively, birth probabilities $p^*_i$), $0 \leq i \leq d - 1$, are positive.

(c) There is a stochastic interpretation of the pgf in \refT{T:Keildisc} even when some of the eigenvalues are negative (see~(4.23) in~\cite{DF}), but we do not know a stochastic proof (i.e., a proof that proceeds by constructing random variables) in that case.

(d) The condition~\eqref{mono+} is closely related to the notion of (stochastic) monotonicity.  All continuous-time, but not all discrete-time, birth-and-death chains are monotone.  In discrete time, monotonicity for a general chain is the requirement that the distributions $P^*(i, \cdot)$ in the successive rows of~$P^*$ be stochastically nondecreasing, i.e.,\ that $\sum_{k > j} P^*(i,k)$ be nondecreasing in~$i$ for each~$j$.  As noted in~\cite{CoxRoesler}, for a discrete-time birth-and-death chain~$P^*$, monotonicity is equivalent to the condition
$$
p^*_{i - 1} + q^*_i \leq 1, \qquad 1 \leq i \leq d.
$$
\end{remark}

We need only prove the discrete-time \refT{T:Keildisc} (or even just the special case where~$P^*$ has positive eigenvalues), for then given a continuous-time birth-and-death generator~$G^*$ we can consider the discrete-time birth-and-death kernels $P^*(\eps) := I + \eps G^*$, where~$I$ denotes the identity matrix and $\eps > 0$ is chosen sufficiently small that $P^*(\eps)$ is nonnegative and has positive eigenvalues.  Let $T(\eps)$ and~$T$ denote the absorption times for $P^*(\eps)$ and~$G^*$, respectively.  Then it is simple to check that $\eps T(\eps)$ converges in law to~$T$; indeed, for any $0 < t < \infty$ we have
$$
\PP(\eps T(\eps) \leq t) = (P^*(\eps))^{\lf t / \eps \rf}(0, d) \to (e^{t G^*})(0, d) = P(T \leq t). 
$$
But the eigenvalues of~$P^*(\eps)$ and of~$-G^*$ are simply related, and suitably scaled geometric random variables converge in law to exponentials, so \refT{T:Keil} follows immediately.

The idea of our proof of \refT{T:Keildisc} is simple:  We show that the absorption time (call it~$T^*$) of~$P^*$ has the same distribution as~$\Th$, where~$\Th$ is the absorption time of a certain pure-birth chain~$\Ph$ whose holding probabilities are precisely the non-unit eigenvalues of~$P^*$. 

We do this by reviewing (in \refS{S:SSD}) and then employing the Diaconis and Fill~\cite{DF} theory of strong stationary duality in discrete time.  In brief, a given absorbing birth-and-death chain~$P^*$ satisfying~\eqref{mono+} is the classical set-valued 
strong stationary dual (SSD) of some monotone birth-and-death chain~$P$ with the same eigenvalues; naturally enough, we will call~$P$ an ``anti-dual'' of~$P^*$.  But, if also the eigenvalues of~$P^*$ are nonnegative, then we show that this~$P$ (and indeed any ergodic birth-and-death chain with nonnegative eigenvalues) in turn also has a pure-birth SSD~$\Ph$ whose holding probabilities are precisely the non-unit eigenvalues of~$P$.  Since we argue that both duals are sharp (i.e.,\ give rise to a stochastically minimal strong stationary time for the $P$-chain), the absorption time~$T^*$ of~$P^*$ has the same distribution as the absorption time~$\Th$ of~$\Ph$, and the latter distribution is manifestly the convolution of geometric distributions.

\begin{remark}
\label{R:path}
(a)~Although our proof of \refT{T:Keildisc} is stochastic, it leaves open [or, rather, left open---see part (c) of this remark] the question of whether the absorption time itself can be represented as an independent sum of explicit geometric random variables; the proof establishes only equality in distribution.  The difficulty with our approach is that there can be many different stochastically minimal strong stationary times for a given chain.

(b)~However, for either of the two steps of our argument we can give sample-path constructions relating the two chains (either~$P^*$ and~$P$, or~$P$ and~$\Ph$).  This has already been carried out in detail for the first step in~\cite{DF}.  For the second step, what this means is that we can show how to watch the $P$-chain $X$ run and contemporaneously construct from it a chain~$\Xh$ with kernel~$\Ph$ in such a way that  the absorption time~$\Th$ of~$\Ph$ is a fastest strong stationary time for~$X$.

(c)~Subsequent to the work leading to the present paper, Diaconis and Miclo~\cite{DM} gave another stochastic proof of~\refT{T:Keil}.  Their proof, which provides an ``intertwining'' between the kernels~$P^*$ and~$\Ph$ (in our notation), yields a construction of exponentials summing to the absorption time, but the construction is, by their own estimation, ``quite involved''.  In a forthcoming paper~\cite{skipfree}, we will exhibit a much simpler such construction, with extensions to skip-free processes.  
\end{remark} 

\refS{S:SSD} is devoted to a brief review of strong stationary duality and a proof that any discrete-time birth-and-death kernel with positive eigenvalues satisfies~\eqref{mono+}.   In \refS{S:anti}, we construct~$P$ from~$P^*$.  In \refS{S:pure} we construct~$\Ph$ from~$P$ and (in \refS{S:path}) describe the sample-path construction discussed in \refR{R:path}(b).  In \refS{S:cont} for completeness we provide continuous-time analogs of our discrete-time auxiliary results, which we find interesting in their own right and which combine to give a direct stochastic proof of \refT{T:Keil}.  \refS{S:RKK} shows how to extend \refT{T:Keil} from the hitting time of state~$d$ to the occupation-time vector for the states $\{0, \dots, d - 1\}$ and connects the present paper with work of Kent~\cite{Kent} and the celebrated Ray--Knight theorem~\cite{Ray, Knight}. 

\section{A quick review of strong stationary duality} \label{S:SSD}

The main purpose of this background section is to review the theory of strong stationary duality only to the extent necessary to understand the proof of~\refT{T:Keildisc}.  For a more general and more detailed treatment, consult~\cite{DF}, especially Sections 2--4.  To a reasonable extent, the notation of this paper matches that of~\cite{DF}.  Strong stationary duality has been used to bound mixing times of Markov chains and also to build perfect simulation algorithms~\cite{interruptible}.

\subsection{Strong stationary duality in general} \label{S:general}

Let~$X$ be an ergodic (irreducible and aperiodic) Markov chain on a finite state space; call its stationary distribution~$\pi$.  A \emph{strong stationary time} is a randomized stopping time~$T$ for~$X$ such that $X_T$ has the distribution~$\pi$ and it independent of~$T$.  Aldous and Diaconis~\cite[Proposition~3.2]{AD87} prove that for any such~$X$ there exists a fastest (i.e.,\ stochastically minimal) strong stationary time, although it is well known that such a fastest time is not (generally) unique.  (Such a fastest time is called a \emph{time to stationarity} in~\cite{DF}, but this terminology has not been widely adopted and so will not be used here.)

A systematic approach to building strong stationary times is provided by the framework of strong stationary duality.  The following specialization of the treatment in Section~2 of~\cite{DF} (see especially Theorem~2.17 and Remark~2.39 there) is sufficient for our purposes.

\begin{theorem}\label{T:SSD}
Let $\pi_0$ and $\pi^*_0$ be probability mass functions on $\{0, 1, \dots d\}$, regarded as row vectors, and let~$P$, $P^*$, and~$\gL$ be transition matrices on~$S$.  Assume that~$P$ is ergodic with stationary distribution~$\pi$, that state~$d$ is absorbing for~$P^*$, and that  the row $\gL(d, \cdot)$ equals~$\pi$.  If $(\pi^*_0, P^*)$ is a \emph{strong stationary dual of $(\pi_0, P)$ with respect to the link~$\gL$} in the sense that
\begin{equation}
\label{algdual}
\pi_0 = \pi^*_0 \gL \quad \mbox{and} \quad \gL P = P^* \gL,
\end{equation}
then there exists a bivariate Markov chain $(X^*, X)$ such that
\begin{enumerate}
\item[{\rm (a)}] $X$ is marginally Markov with initial distribution~$\pi_0$ and transition matrix~$P$;
\item[{\rm (b)}] $X^*$ is marginally Markov with initial distribution~$\pi^*_0$ and transition matrix~$P^*$;
\item[{\rm (c)}] the absorption time~$T^*$ of~$X^*$ is a strong stationary time for~$X$.
\end{enumerate} 
Moreover, if $\gL(i, d) = 0$ for $i = 0, \dots, d - 1$, then the dual is \emph{sharp} in the sense that~$T^*$ is a fastest strong stationary time for~$X$.  
\end{theorem}

\begin{remark}
\label{R:special}
In both our applications of \refT{T:SSD} (Sections~\ref{S:anti} and~\ref{S:pure}),
\begin{enumerate}
\item the initial distributions $\pi_0$ and $\pi^*_0$ are both taken to be unit mass $\delta_0$ at~$0$, and $\gL(0, \cdot) = \delta_0$, too, so only the second equation in~\eqref{algdual} needs to be checked; and
\item the link~$\gL$ is lower triangular, from which we observe that the corresponding dual is sharp and (if also the diagonal elements of~$\gL$ are all positive) that, given~$P$, there is at most one stochastic matrix~$P^*$ satisfying~\eqref{algdual}, namely, $P^* = \gL P \gL^{-1}$.
\end{enumerate}
\end{remark}

\subsection{Classical (set-valued) strong stationary duals} \label{S:set}

Let~$P$ be ergodic with stationary distribution~$\pi$, and let~$H$ denote the corresponding cumulative distribution function (cdf):
$$
H_j = \sum_{i \leq j} \pi_i.
$$
Let~$\gL$ be the link of truncated stationary distributions:
\begin{equation}
\label{trunc}
\gL(x^*, x) = {\bf 1}(x \leq x^*) \pi_x / H_{x^*}.
\end{equation}
If~$P$ is a monotone birth-and-death chain (more generally, if~$P$ is arbitrary and the time reversal of~$P$ is monotone---see~\cite[Theorem 4.6]{DF}), then a dual~$P^*$ exists [and is sharp and unique by \refR{R:special}(ii)]:

\begin{theorem}
\label{T:classic}
Let~$P$ be a monotone ergodic birth-and-death chain on $\{0, \dots, d\}$ with stationary cdf~$H$.  Then~$P$ has a sharp strong stationary dual~$P^*$ with respect to the link of truncated stationary distributions.  The chain~$P^*$ is also birth-and-death, with death, hold, and birth probabilities (respectively)
\begin{equation}
\label{BDdual}
q^*_i = \frac{H_{i - 1}}{H_i} p_i  \qquad r^*_i = 1 - (p_i + q_{i + 1}), \qquad p^*_i = \frac{H_{i + 1}}{H_i} q_{i + 1} .
\end{equation}
\end{theorem}

See Sections~3--4 of~\cite{DF} for an explanation as to why the dual in~\refT{T:classic} is called ``set-valued''; in this paper we shall refer to it as the ``classical'' SSD.  The equations~\eqref{BDdual} reproduce~\cite[(4.18)]{DF}.

\subsection{Positivity of eigenvalues and stochastic monotonicity for birth-and-death chains}
\label{S:mono}

When we prove \refT{T:Keildisc} assuming that~$P^*$ has positive eigenvalues, we will utilize the strengthened monotonicity condition~\eqref{mono+}.  Part~(a) of the following lemma provides justification.

\begin{lemma} \label{L:mono}
Let~$P^*$ be the kernel of any birth-and-death chain on $\{0, \dots, d\}$.
\begin{enumerate}
\item[{\rm (a)}] If~$P^*$ has positive eigenvalues, then~\eqref{mono+} holds.
\item[{\rm (b)}] If~$P^*$ has nonnegative eigenvalues, then~$P^*$ is monotone.
\end{enumerate}
\end{lemma}

\begin{proof}
(a)~By perturbing~$P^*$ if necessary, we may assume that~$P^*$ is ergodic.  Then~$P^*$ is diagonally similar to a positive definite matrix whose principal minor corresponding to rows and columns $i - 1$ and~$i$ is $r^*_{i - 1} r^*_i - p^*_{i - 1} q^*_i$, so 
$$
0 < r^*_{i - 1} r^*_i - p^*_{i - 1} q^*_i \leq (1 - p^*_{i - 1})(1 - q^*_i) - p^*_{i - 1} q^*_i = 1 - p^*_{i - 1} - q^*_i.
$$

(b)~This follows by perturbation from part~(a).
\end{proof}

\begin{remark} \label{noconverse}
Both converse statements are false.  For any given $d \geq 2$, the condition~\eqref{mono+} does not imply nonnegativity of eigenvalues, not even for chains~$P^*$ satisfying the hypotheses of \refT{T:Keildisc}.  An explicit counterexample for $d = 2$ is
$$
P^* = 
\left[ \begin{array}{lll}
           0.50 & 0.50 & 0 \\
           0.49 & 0.02 & 0.49 \\
           0       & 0      & 1
         \end{array}                     \right],
$$
whose smallest eigenvalue is $(26 - \sqrt{3026}) / 100 < 0$.  For general $d \geq 2$, perturb the direct sum of this counterexample with the identity matrix. 
\end{remark}         
  
\section{An anti-dual~$P$ of the given~$P^*$} \label{S:anti}

As discussed in~\refS{S:intro}, the main discrete-time theorem, \refT{T:Keildisc}, follows from the chief results, Theorems~\ref{T:anti} and~\ref{T:spectraldual}, of this section and the next. 

Under the strengthened monotonicity condition~\eqref{mono+} (with no assumption here about nonnegativity of the eigenvalues), the anti-dual construction of~\refT{T:anti} exhibits the given chain (call its kernel~$P^*$) as the classical SSD of another birth-and-death chain.

\begin{theorem}\label{T:anti}
Consider a discrete-time birth-and-death chain $P^*$ on $\{0, \dots, d\}$ started at~$0$, and suppose that~$d$ is an absorbing state.  Write $q^*_i$, $r^*_i$, and $p^*_i$ for its death, hold, and birth probabilities, respectively.  Suppose that $p^*_i > 0$ for $0 \leq i \leq d - 1$, that $q^*_i > 0$ for $1 \leq i \leq d-1$, and that $p^*_{i - 1} + q^*_i < 1$ for $1 \leq i \leq d$.  Then $P^*$ is the classical (and hence sharp) SSD of some monotone ergodic birth-and-death kernel~$P$ on $\{0, \dots, d\}$. 
\end{theorem}

\begin{proof}
\noqed
In light of \refR{R:special}(i), we have dispensed with initial distributions.  The claim is that $P^*$ is related to some monotone ergodic~$P$ with stationary cdf~$H$ via~\eqref{BDdual}.  We will \emph{begin} our proof by defining a suitable function~$H$, and then we will construct~$P$.

We inductively define a strictly increasing function $H: \{0, \dots, d\} \to (0, 1]$.  Let $H_d :=1$, and define $H_{d - 1} \in (0, 1)$ in (for now) arbitrary fashion.  Having defined $H_d, \dots, H_i$ (for some $1 \leq i \leq d - 1$), choose the value of $H_{i - 1} \in (0, H_i)$ so that
\begin{equation} 
\label{Hdef} 
\left( \frac{H_i}{H_{i - 1}} - 1 \right) q^*_i = \left( 1 - \frac{H_i}{H_{i + 1}} \right) p^*_i;
\end{equation}
this is clearly possible since the right side of~\eqref{Hdef} is in $(0, 1)$ and the left side, as a function of the variable $H_{i - 1}$, decreases from $\infty$ at $H_{i - 1} = 0+$ to~$0$ at $H_{i - 1} = H_i$.  It is also clear that by choosing $H_{d - 1}$ sufficiently close to~$1$, we can make \emph{all} the ratios $H_i / H_{i - 1}$ ($i = 1, \dots d$) as (uniformly) close to~$1$ as we wish.

Next, define $q_0 := 0$,
\begin{equation}
\label{p0def}
p_0 := \left( 1 - \frac{H_0}{H_1} \right) p^*_0,
\end{equation}
and, for $1 \leq i \leq d$,
\begin{equation}
\label{pqdef}
p_i := \frac{H_i}{H_{i - 1}} q^*_i,  \qquad q_i := \frac{H_{i - 1}}{H_i} p^*_{i - 1}.
\end{equation}
When the $H$-ratios are taken close enough to~$1$, then for $0 \leq i \leq d$ we have $p_i + q_i < 1$ and we define
$$
r_i := 1 - p_i - q_i > 0.
$$
The kernel~$P$ with death, hold, and birth probabilities $q_i$, $r_i$, and $p_i$ is irreducible and aperiodic, and thus ergodic.  To complete the proof, will also show
\begin{enumerate}
\item[(a)] $P$ is monotone (recall:\ equivalent to $p_i + q_{i + 1} \leq 1$ for $0 \leq i \leq d - 1$),
\item[(b)] $P$ has stationary cdf~$H$, and
\item[(c)] $P^*$ is the classical SSD of~$P$.
\end{enumerate}

For~(a) we simply observe, using~\eqref{pqdef} and~\eqref{Hdef}, that 
\begin{equation}
\label{a}
p_i + q_{i + 1} = \frac{H_i}{H_{i - 1}} q^*_i + \frac{H_i}{H_{i + 1}} p^*_i = q^*_i + p^*_i \leq 1
\end{equation}
for $1 \leq i \leq d - 1$; and similarly that
$$
p_0 + q_1 =  \left( 1 - \frac{H_0}{H_1} \right) p^*_0 + \frac{H_0}{H_1} p^*_0 = p^*_0 \leq 1.
$$

For~(b) we observe, again using~\eqref{pqdef} and~\eqref{Hdef}, that the detailed balance condition
$$
(H_i - H_{i - 1}) p_i = (H_i - H_{i - 1}) \frac{H_i}{H_{i - 1}} q^*_i = (H_{i + 1} - H_i) \frac{H_i}{H_{i + 1}} p^*_i = (H_{i + 1} - H_i) q_{i + 1}
$$
holds for $1 \leq i \leq d - 1$; by~\eqref{p0def} and~\eqref{pqdef}, it also holds for $i = 0$:
$$
H_0 p_0 = (H_1 - H_0) \frac{H_0}{H_1} p^*_0 = (H_1 - H_0) q_1.
$$

For~(c), we simply verify that~\eqref{BDdual} holds:\ for $0 \leq i \leq d$ (with $H_{-1} := 0$), from~\eqref{pqdef} and~\eqref{a},
\begin{equation}
\label{dual}
\frac{H_{i - 1}}{H_i} p_i = q^*_i, \qquad \frac{H_{i + 1}}{H_i} q_{i + 1} = p^*_i, \qquad p_i + q_{i + 1} = q^*_i + p^*_i = 1 - r^*_i. \qquad ~\qedsymbol
\end{equation}
\end{proof}

\begin{remark}
\label{R:antirem}
Once the value of $H_{d - 1}$ is chosen, the definitions of~$H$ and~$P$ are forced; indeed, if the detailed balance condition and~\eqref{dual} are to hold, then we must have~\eqref{Hdef}--\eqref{pqdef}.
\end{remark}

\section{A pure birth ``spectral'' dual of~$P$} \label{S:pure}

In this section we construct a sharp pure birth dual~$\Ph$ for any ergodic birth-and-death chain~$P$ on $\{0, \dots d\}$ with nonnegative eigenvalues started in state~$0$.  
When this construction is applied in the proof of~\refT{T:Keildisc} to the chain~$P$ resulting from~$P^*$ by application of \refT{T:anti}, assuming nonnegativity of the eigenvalues of~$P^*$ yields the required nonnegativity of the eigenvalues of~$P$ in \refT{T:spectraldual}; indeed, as noted in \refR{R:special}(ii), the matrices~$P$ and~$P^*$ are similar.
\ignore{IGNORE THIS FROM EARLIER DRAFT BECAUSE IT HAS BEEN IMPROVED!
Note that it is guaranteed that the chain~$P$ constructed in the proof of \refT{T:anti} will have nonnegative eigenvalues, provided the $H$-ratios are taken close enough to~$1$, when
$$
q^*_i + p^*_{i - 1} < 1 / 2,
$$
for all~$i$, since then $r_i \geq 1/2$ for all~$i$.}

Our construction of the pure birth dual specializes a SSD construction of Matthews \cite{Matt} for general reversible chains with nonnegative eigenvalues; that construction is closely related to the spectral decomposition of the transition matrix.  For completeness and the reader's convenience, and because for birth-and-death chains (a)~we can give a more streamlined presentation with minimal reference to eigenvectors and (b)~we wish to establish the new result that the resulting dual is sharp, we do not presume familiarity with~\cite{Matt}.

To set up our construction we need some notation.  Let~$P$ be an ergodic birth-and-death chain on $\{0, \dots, d\}$ with stationary probability mass function~$\pi$ (note that~$\pi$ is everywhere positive) 
and nonnegative eigenvalues, say $0 \leq \theta_0 \leq \theta_1 \leq \dots \leq \theta_{d - 1} < \theta_d = 1$.  (It is well known~\cite{Keil}~\cite[Theorem~4.20]{DF} that the eigenvalues are all distinct, but we will not need this fact.)  Let~$I$ denote the identity matrix and define
\begin{equation}
\label{Qkdef}
Q_k := (1 - \theta_0)^{-1} \cdots (1 - \theta_{k-1})^{-1} (P - \theta_0 I) \cdots (P - \theta_{k-1} I), \quad k = 0, \dots, d,
\end{equation}
with the natural convention $Q_0 := I$.  Note that for $k = 0, \dots, d - 1$ we have
\begin{equation}
\label{Qevo}
Q_k P = \theta_k Q_k + (1 - \theta_k) Q_{k + 1}.
\end{equation}

\begin{lemma}
\label{L:spectral}
The matrices $Q_k$ are all stochastic, and every row of $Q_d$ equals~$\pi$.
\end{lemma}

\begin{proof}
For the first assertion it is clear that the rows of~$Q_k$ all sum to~$1$, so the only question is whether $Q_k$ is nonnegative.  But $P = D^{- 1 / 2} S D^{1 / 2}$, where $D = \diag(\pi)$ and~$S$ is symmetric, so the nonnegativity of $Q_k$ follows from that of 
$$
S_k := (S - \theta_0 I) \cdots (S - \theta_{k - 1} I),
$$ 
which in turn is an immediate consequence of (the rather nontrivial) Theorem~3.2 in~\cite{MW} using only that~$S$ is nonnegative and symmetric.

For the second assertion, write
$$
S = \sum_{r = 0}^d \theta_r u_r u_r^T,
$$
where the column vectors $u_0, \dots, u_d$ form an orthogonal matrix and $u_d$ has $i$th entry $\sqrt{\pi_i}$.  Then, as noted at~(2.6) of~\cite{MW},
$$
S_k = \sum_{r = k}^d \left[ \prod_{t = 0}^{k - 1} (\theta_r - \theta_t) \right] u_r u_r^T.
$$
In particular, $S_d = (1 - \theta_0) \cdots (1 - \theta_{d - 1}) u_d u_d^T$, so every row of $Q_d$ equals~$\pi$.
\end{proof}

Now let $\delta_0$ denote unit mass at~$0$ (regarded as a row vector), and define the probability mass functions
\begin{equation}
\label{glkdef}
\gl_k := \delta_0 Q_k, \quad k = 0, \dots, d.
\end{equation}
Let $\gLh$ [so named to distinguish it from the classic link~$\gL$ of~\eqref{trunc}] be the lower-triangular square matrix with successive rows $\gl_0, \dots, \gl_d$, and define~$\Ph$ to be the pure-birth chain transition matrix on $\{0, \dots, d\}$ with holding probability $\theta_i$ at state~$i$ for $i = 0, \dots, d$; that is,
\begin{equation}
\label{phat}
\ph_{i j} := 
 \begin{cases}
 \theta_i   & \mbox{if $j = i$} \\
 1 - \theta_i   & \mbox{if $j = i + 1$} \\
 0   & \mbox{otherwise}.
 \end{cases} 
\end{equation}

\begin{theorem}
\label{T:spectraldual}
Let~$P$ be an ergodic birth-and-death chain on $\{0, \dots, d\}$ with nonnegative eigenvalues.  In the above notation, $\Ph$ is a sharp strong stationary dual of~$P$ with respect to the link~$\gLh$.
\end{theorem} 

\begin{proof}
We have again dispensed with initial distributions by \refR{R:special}(i).  The desired equation $\widehat{\gL} P = \Ph \gLh$ is equivalent to
$$
\gl_k P = \theta_k \gl_k + (1 - \theta_k) \gl_{k + 1}, \quad k = 0, \dots, d - 1; \qquad \gl_d P = \gl_d,
$$
which is true because $\gl_d = \pi$ and, for $k = 0, \dots, d - 1$,
$$
\gl_k P = \delta_0 Q_k P = \theta_k \gl_k + (1 - \theta_k) \gl_{k + 1}
$$
by~\eqref{Qevo}.  The SSD is sharp because~$\gLh$ is lower triangular; recall \refR{R:special}(ii). 
\end{proof}

\begin{remark}
\refL{L:spectral} is interesting and, as we have now seen, gives rise to the construction of a new ``spectral'' SSD for a certain subclass of monotone birth-and-death chains, namely, chains with nonnegative eigenvalues [recall \refL{L:mono}(b)].  But for the proof of \refT{T:Keildisc} one could make do without the nonnegativity of the matrix~$\gLh$ by taking the approach of Matthews~\cite{Matt} and considering the chain~$P$ started in a suitable mixture of $\delta_0$ and the stationary distribution~$\pi$.  We omit further details.
\end{remark}

\subsection{Sample-path construction of the spectral dual} \label{S:path}

Let~$X$ be an ergodic birth-and-death chain on $\{0, \dots, d\}$ with kernel~$P$ having nonnegative eigenvalues, assume $X_0 = 0$, and let~$T$ be any fastest strong stationary time for~$X$.  Independent of interest in Theorems~\ref{T:Keil} and \ref{T:Keildisc}, \refT{T:spectraldual} gives the first stochastic interpretation of the individual geometrics in the representation of the distribution of~$T$ as a convolution of geometric distributions.  In this subsection we carry this result one step further by showing how to construct, sample path by sample path, a particular fastest strong stationary time~$\Th$ which is the sum of explicitly identified independent geometric random variables.

The idea is simple.  \refT{T:spectraldual} shows that~$\Ph$ of~\eqref{phat} is an ``algebraic'' dual of~$P$ in the sense that the matrix-equation $\gLh P = \Ph \gLh$ holds.  But whenever algebraic duality holds for any finite-state ergodic chain with respect to any link ($\gLh$ in our case), Section~2.4 of~\cite{DF} shows explicitly how to construct, from~$X$ and independent randomness, a dual Markov chain ($\Xh$ in our case, with kernel~$\Ph$) such that the absorption time~$\Th$ of~$\Xh$ is a strong stationary time for~$X$; since~$\gLh$ is lower triangular, $\Th$ will be stochastically optimal.  So to describe our construction of~$\Xh$ (and hence~$\Th$) we need only specialize the construction of \cite[Section~2.4]{DF} [see especially~(2.36) there]. 

The chain~$X$ starts with $X_0 = 0$ and we set $\Xh_0 = 0$.  Inductively, we will have $\gLh(\Xh_t, X_t) > 0$ (and so $X_t \leq \Xh_t$) at all times~$t$.  The value we construct for $\Xh_t$ depends only on the values of~$\Xh_{t - 1}$ and~$X_t$ and independent randomness.  Indeed, given $\Xh_{t - 1} = \xh$ and $X_t = y$, if $y \leq \xh$ then our construction sets $\Xh_t = \xh + 1$ with probability
\begin{equation}
\label{spectralpath}
\frac{\Ph(\xh, \xh + 1) \gLh(\xh + 1, y)}{(\Ph \gLh)(\xh, y)} = \frac{(1 - \theta_{\xh}) \gLh(\xh + 1, y)}{\theta_{\xh} \gLh(\xh, y) + (1 - \theta_{\xh}) \gLh(\xh + 1, y)} =  \frac{(1 - \theta_{\xh}) Q_{\xh + 1}(0, y)}{(Q_{\xh} P)(0, y)}
\end{equation}
and $\Xh_t = \xh$ with the complementary probability; if $y = \xh + 1$ (which is the only other possibility, since $y = X_t \leq X_{t - 1} + 1 \leq \xh + 1$ by induction), then we set $\Xh_t = \xh + 1$ with certainty.

The independent geometric random variables, with sum~$\Th$, are the waiting times between successive births in the chain~$\Xh$ we have built.  Thus it is no longer true that the individual geometric distributions ``have no known interpretation in terms of the underlying [ergodic] birth and death chain''  \cite[Section 4, Remark 1]{DS}; likewise, for continuous time consult \refS{S:pathc} herein.

\begin{example}
{\rm 
Consider the well-studied Ehrenfest chain, with holding probability $1 / 2$:
$$
q_i = \frac{i}{2 d}, \quad r_i = \frac{1}{2}, \quad p_i = \frac{d - i}{2 d}, \qquad i = 0, \dots, d.
$$
The eigenvalues are $\theta_i \equiv i / d$.  A straightforward proof by induction using~\eqref{glkdef} and~\eqref{Qevo} confirms that $\gl_k$ is the binomial distribution with parameters~$k$ and~$1 / 2$:
\begin{equation}
\label{binomial}
\gLh(\xh, x) \equiv {\xh \choose x} 2^{- \xh}.
\end{equation}
Thus the probability~\eqref{spectralpath} reduces to
$$
\frac{(d - \xh) (\xh + 1)}{2 \xh (\xh + 1 - y) + (d - \xh)(\xh + 1)}.
$$

The chain we have described lifts naturally to random walk on the set ${\bf Z}^d_2$ of binary $d$-tuples whereby one of the~$d$ coordinates is chosen uniformly at random and its entry is then replaced randomly by~$0$ or~$1$.  It is interesting to note that the sharp pure-birth SSD chain constructed in this example does \emph{not} correspond to the well-known ``coordinate-checking'' sharp SSD (see Example~3.2 of~\cite{DF}).  Indeed, expressed in the birth-and-death chain domain, the coordinate-checking dual is a pure-birth chain, call it $\Xh'$, such that the construction of~$\Xh'_t$ depends not only on~$\Xh'_{t - 1}$ and~$X_t$ but also on~$X_{t - 1}$.  The construction rules are that if $\Xh'_{t - 1} = \xh$, $X_{t - 1} = x$, and $X_t = y$, then~$\Xh'_t$ is set to $\xh + 1$ with probability
$$
\mbox{$0$ if $y = x - 1$,\ \ \ \ $1 - \frac{\xh}{d}$ if $y = x$,\ \ \ \ and\ \ \ \ $\frac{d - \xh}{d - x}$ if $y = x + 1$,}
$$
and otherwise~$\Xh'_t$ holds at~$\xh$.  Both duals correspond to the same link~\eqref{binomial} and the (marginal) transition kernels for~$\Xh$ and~$\Xh'$ are the same, but the bivariate constructions of $(\Xh, X)$ and $(\Xh', X)$ are different.  

The freedom for such differences was noted in \cite[Remark 2.23(c)]{DF} and exploited in the creation of an interruptible perfect simulation algorithm (see \cite[Remark~9.8]{interruptible}).  In fact, $\Xh'$ (when lifted to ${\bf Z}^d_2$) corresponds to the construction used in~\cite{interruptible}.  An advantage of the $\Xh$-construction of the present paper is that it allows (both in our Ehrenfest example and in general) for holding probabilities that are arbitrary (subject to nonnegativity of eigenvalues); in the paragraph containing~\eqref{spectralpath}, all that changes when a weighted average of the transition kernel and the identity matrix is taken are the eigenvalues $\theta_0, \dots, \theta_{d - 1}$.
}
\end{example}

\section{Continuous-time analogs of other results} \label{S:cont}

As discussed in \refS{S:intro}, the continuous-time \refT{T:Keil} follows immediately from the discrete-time \refT{T:Keildisc}.  Another way to prove \refT{T:Keil} is to repeat the proof of \refT{T:Keildisc} by establishing continuous-time analogs (namely, the next three results) of the auxiliary results (\refT{T:anti}, \refL{L:spectral}, and \refT{T:spectraldual}) in the preceding two sections; we find these interesting in their own right.  The continuous-time results are easy to prove utilizing the continuous-time SSD theory of~\cite{dualityc}, either by repeating the discrete-time proofs or by applying the discrete-time results to the appropriate kernel $P^*(\eps) = I + \eps G^*$ or $P(\eps) = I + \eps G$, with $\eps > 0$ chosen sufficiently small to meet the hypotheses of those results; so we state the results without proof.

In \refS{S:pathc} we will present the analog of \refS{S:path} for continuous time.

Here, first, is the analog of \refT{T:anti}.

\begin{theorem}\label{T:antic}
Consider a continuous-time birth-and-death chain with generator~$G^*$ on $\{0, \dots, d\}$ started at~$0$, and suppose that~$d$ is an absorbing state.  Write $\mu^*_i$ and $\gl^*_i$ for its death and birth rates, respectively.  Suppose that $\gl^*_i > 0$ for $0 \leq i \leq d - 1$ and that $\mu^*_i > 0$ for $1 \leq i \leq d-1$.  Then $G^*$ is the classical set-valued (and hence sharp) SSD of some ergodic birth-and-death generator~$G$ on $\{0, \dots, d\}$. 
\end{theorem}
\ignore{
That is, the claim is that~$G^*$ is related to some ergodic birth-and-death generator~$G$ via (3.12) of \cite{dualityc}.  There, $H$ denotes the stationary cdf for the birth-and-death chain with generator~$G$.
Define~$H$ as in the proof of the discrete-time \refT{T:anti}, but with $(\mu^*_i, \gl^*_i)$ playing the roles of $(q^*_i, p^*_i)$ in~\eqref{Hdef}.  Next, in analogous fashion as for discrete time, define $\gl_0 := 0$,
\begin{equation}
\label{gl0def}
\gl_0 := \left( 1 - \frac{H_0}{H_1} \right) \gl^*_0,
\end{equation}
and, for $1 \leq i \leq d$,
\begin{equation}
\label{pqdef}
\gl_i := \frac{H_i}{H_{i - 1}} \mu^*_i,  \qquad \mu_i := \frac{H_{i - 1}}{H_i} \gl^*_{i - 1}.
\end{equation}
We have thus defined the generator~$G$ of a birth-and-death chain with death and birth rates $\mu_i$ and $\gl_i$, respectively.  Just as in discrete time, $G$ has stationary cdf~$H$, and~$G^*$ is the usual (classical set-valued) SSD of~$G$.
}

To set up the second result we need a little notation.  Let~$G$ be the generator of a continuous-time ergodic birth-and-death chain on $\{0, \dots, d\}$ with stationary probability mass function~$\pi$ and eigenvalues $\nu_0\geq \nu_1 \geq \dots \geq \nu_{d - 1} > \nu_d = 0$ for $- G$.  (Again, we don't need the fact~\cite{Keil} that the eigenvalues are distinct.)  Define
\begin{equation}
\label{Qkdefc}
Q_k := \nu_0^{-1} \cdots \nu_{k - 1}^{-1} (G + \nu_0 I) \cdots (G + \nu_{k-1} I), \quad k = 0, \dots, d,
\end{equation}
with the natural convention $Q_0 := I$.
\begin{lemma}
\label{L:spectralc}
The matrices $Q_k$ are all stochastic, and every row of $Q_d$ equals~$\pi$.
\end{lemma}

Now define~$\widehat{\gL}$ in terms of the $Q_k$'s as in the paragraph preceding \refT{T:spectraldual}, and let~$\Gh$ be the pure-birth generator on $\{0, \dots, d\}$ with birth rate $\nu_i$ at state~$i$ for $i = 0, \dots, d$.  

\begin{theorem}
\label{T:spectraldualc}
Let~$G$ be the generator of an ergodic birth-and-death chain on $\{0, \dots, d\}$.  In the above notation, $\Gh$ is a sharp strong stationary dual of~$G$ with respect to the link~$\widehat{\gL}$:
$$
\widehat{\gL} G = \Gh \widehat{\gL}.
$$
\end{theorem}

\subsection{Sample-path construction of the continuous-time spectral dual} \label{S:pathc}

Let~$X$ be an ergodic continuous-time birth-and-death chain on $\{0, \dots, d\}$, adopt all the notation of \refS{S:cont} thus far, and assume $X_0 = 0$.   In this subsection by a routine application of Section~2.3 of~\cite{dualityc} we give a simple sample-path construction of a ``spectral dual'' pure birth chain~$\Xh$ with generator~$\Gh$ as described just before \refT{T:spectraldualc}; its absorption time~$\Th$ is then a fastest strong stationary time for~$X$ and the independent exponential random variables with sum~$\Th$ are simply the waiting times for the successive births for~$\Xh$.  We thus obtain a stochastic proof, with explicit identification of individual exponential random variables, of Theorem~5 in~\cite{dualityc}.

The chain~$X$ starts with $X(0) = 0$ and we set $\Xh(0) = 0$.  Let $n \geq 1$ and suppose that~$\Xh$ has been constructed up through the epoch $\tau_{n - 1}$ of the $(n - 1)$st transition for the bivariate process $(\Xh, X)$; here $\tau_0 := 0$.  We describe next, in terms of an exponential random variable $\Vh$, how to define~$\tau_n$ and $\Xh(\tau_n)$; we will have $\gLh(\Xh(\tau_n), X(\tau_n)) > 0$ and hence $X(\tau_n) \leq \Xh(\tau_n)$.  Write $(\xh, x)$ for the value of $(\Xh, X)$ at time $\tau_{n - 1}$; by induction we have $\gLh(\xh, x) > 0$.

Let $\Vh_n$ be exponentially distributed with rate
\begin{equation}
\label{rate}
r = \nu_{\xh} \gLh(\xh + 1, x) / \gLh(\xh, x),
\end{equation}
independent of $\Vh_1, \dots, \Vh_{n - 1}$ and the chain~$X$.  Consider two (independent) exponential waiting times begun at epoch~$\tau_{n - 1}$:\ a first for the next transition of the chain~$X$, and a second with rate~$r$.  How we proceed breaks into two cases: 
\begin{enumerate}
\item If the first waiting time is smaller than the second, then~$\tau_n$ is the epoch of this next transition for~$X$ and we set $\Xh(\tau_n) = \xh = \Xh(\tau_{n - 1})$ (with certainty) except in one circumstance:\ if $X(\tau_n) = \xh + 1$, then we set $\Xh(\tau_n) = \xh + 1$, too.
\item If the second waiting time is smaller, then $\tau_n = \tau_{n - 1} + \Vh_n$ and we set $\Xh(\tau_n) = \xh + 1$.
\end{enumerate}

\begin{example}
{\rm 
Consider the continuous-time version of the Ehrenfest chain with death rates $\mu_i \equiv i$ and birth rates $\gl_i \equiv d - i$, $0 \leq i \leq d$; the eigenvalues are $\nu_i \equiv 2 (d - i)$.  Then~$\gLh$ is again the link~\eqref{binomial} of binomial distributions, and the rate~\eqref{rate} reduces to
$$
r = \frac{(d - \xh) (\xh + 1)}{\xh + 1 - x}.
$$ 
}
\end{example}

\section{Occupation times and connection with Ray--Knight Theorem}
\label{S:RKK}

Our final section utilizes work of Kent~\cite{Kent}; see the historical note at the end of Section~1 of \cite{DM} for closely related material.  We show how to extend the continuous-time \refT{T:Keil} from the hitting time of state~$d$ first to the 
occupation-time vector for the states $\{0, \dots, d - 1\}$ and then to the the local time of Brownian motion, thereby proving the Ray--Knight theorem~\cite{Ray, Knight}. 

\subsection{From hitting time to occupation times}

Consider a continuous-time irreducible birth-and-death chain with generator $G^*$.  It is then immediate from the Karlin--McGregor theorem (\refT{T:Keil}) that the hitting time $T^*$ of 
state~$d$ has Laplace transform
\begin{equation}
\label{htLT}
{\bf E}\,e^{-u T^*} = \frac{\det(- G_0)}{\det(- G_0 + u I)},
\end{equation}
with $G_0$ obtained from $G^*$ by leaving off the last row and column.

Equation~\eqref{htLT} gives the distribution of the total time elapsed before the chain hits state~$d$.  But how is that time apportioned to the states $0, \dots, d - 1$?  This question can be answered from~\eqref{htLT} using a neat trick of Kent~\cite{Kent} [see the last sentence of his Remark~(1) on page~164].  To find the multivariate distribution of the occupation-time vector ${\bf T} = (T_0, T_1, \dots, T_{d - 1})$, where $T_i$ denotes the occupation time of (i.e.,\ amount of time spent in) state~$i$, it of course suffices to compute the value ${\bf E}\,e^{- \langle {\bf u}, {\bf T} \rangle}$ of the Laplace transform for any vector ${\bf u} = (u_0, \dots, u_{d - 1})$ with strictly positive entries.  But the distribution of the random variable $\langle {\bf u}, {\bf T} \rangle = \sum u_i T_i$ is that of the time to absorption for the time-changed generator $G_{\bf u}$ (say) obtained by dividing the $i$th row of $G^*$ by $u_i$ for $i = 0, \dots, d - 1$.  Therefore, by~\eqref{htLT} and the scaling property of determinants,
$$
 {\bf E}\,e^{- \langle {\bf u}, {\bf T} \rangle} = \frac{\det(- G_{\bf u})}{\det(- G_{\bf u} + I)} = \frac{\det(- G_0)}{\det(- G_0 + U)},
$$
where $U := \mbox{diag}(u_0, \dots, u_{d - 1})$.

\subsection{From occupation times to the Ray--Knight Theorem}

Call the stationary distribution~$\pi$.  Then the matrix $S := D(- G_0) D^{-1}$ is (strictly) positive definite, where $D := \mbox{diag}(\sqrt{\pi})$.  Let  $\Sigma := \frac{1}{2} S^{-1}$.  By direct calculation, ${\bf T}$ has the same law as ${\bf Y} + {\bf Z}$, where ${\bf Y}$ and ${\bf Z}$ are independent random vectors with the same law and ${\bf Y}$ is the coordinate-wise square of a Gaussian random vector ${\bf V} \sim {\rm N}(0, \Sigma)$.

Kent~\cite{Kent} uses and extends this ``double derivation'' of ${\cal L}({\bf T})$ to prove the theorem of Ray~\cite{Ray} and Knight~\cite{Knight} expressing the local time of Brownian motion as the sum of two independent 2-dimensional Bessel processes (i.e., as the sum of two independent squared Brownian motions).
\bigskip  

{\bf Acknowledgments.\ }We thank Persi Diaconis for helpful discussions, and Raymond Nung-Sing Sze and Chi-Kwong Li for pointing out the reference~\cite{MW}.


\begin{thebibliography}{99}
\def\nobibitem#1\par{}

\bibitem{AD87}
Aldous, D.\ and Diaconis, P.\ \ 
Strong uniform times and finite random walks.
\emph{Adv.\ in Appl.\ Math.} {\bf 8} (1987), 69--97.

\bibitem{BS}
Brown, M.\ and Shao, Y.~S.\ \ 
Identifying coefficients in the spectral representation for first passage time distributions.
\emph{Probab.\ Eng.\ Inform.\ Sci.} {\bf 1} (1987), 69--74.

\bibitem{CoxRoesler}
Cox, J.~T.\ and \Roesler, U.\ \ 
A duality relation for entrance and exit laws for Markov monotone Markov processes.
\emph{Ann.\ Probab.} {\bf 13} (1985), 558--565.

\bibitem{DF}
Diaconis, P.\ and Fill, J.\ \ 
Strong stationary times via a new form of duality.
\emph{Ann.\ Probab.} {\bf 18} (1990), 1483--1522.

\bibitem{DM}
Diaconis, P.\ and Miclo, L.\ \ 
On times to quasi-stationarity for birth and death processes.
\emph{J.\ Theoret.\ Probab.} (2009), to appear.

\bibitem{DS}
Diaconis,~P.\ and Saloff-Coste,~L.\ \ 
Separation cut-offs for birth and death chains.
\emph{Ann.\ Appl.\ Probab.} {\bf 16} (2006), 2098--2122.

\bibitem{dualityc}
Fill, J.~A.\ \ 
Strong stationary duality for continuous-time Markov chains.  Part~I:\ Theory.
\emph{J.\ Theoret.\ Probab.} {\bf 5} (1992), 45--70.

\bibitem{interruptible}
Fill, J.~A.\ \ 
An interruptible algorithm for perfect sampling via Markov chains.
\emph{Ann.\ Appl.\ Probab.} {\bf 8} (1998), 131--162.

\bibitem{skipfree}
Fill, J.~A.\ \ 
On hitting times and fastest strong stationary times for skip-free and more general chains.
\emph{J.\ Theoret.\ Probab.} (2009), to appear.

\bibitem{KM}
Karlin, S.\ and McGregor, J.\ \ 
Coincidence properties of birth and death processes.
\emph{Pacific J.\ Math.} {\bf 9} (1959), 1109--1140.

\bibitem{Keillog}
Keilson, J.\ \ 
Log-concavity and log-convexity in passage time densities for of diffusion and birth-death processes.
\emph{J.\ Appl.\ Probab.} {\bf 8} (1971), 391--398.

\bibitem{Keil}
Keilson, J.\ \ 
\emph{Markov Chain Models---Rarity and Exponentiality.}
Springer, New York, 1979.

\bibitem{Kent}
Kent, J.~T.\ \ 
The appearance of a multivariate exponential distribution in sojourn times for birth-death and diffusion processes. 
In \emph{Probability, statistics and analysis}, volume~79 of
\emph{London Math.\ Soc.\ Lecture Note Ser.},\ pages 161--179, Cambridge Univ.\ Press, Cambridge, 1983.

\bibitem{Knight}
Knight, F.~B.\ \ 
Random walks and a sojourn density process of Brownian motion. 
\emph{Trans.\ Amer.\ Math.\ Soc.} {\bf 109} (1963), 56--86. 

\bibitem{Matt}
Matthews, P.\ \
Strong stationary times and eigenvalues.
\emph{J.\ Appl.\ Probab.} {\bf 29} (1992), 228--233.

\bibitem{MW}
Micchelli, C.~A.\ and Willoughby, R.~A.\ \ 
On functions which preserve the class of Stieltjes matrices.
\emph{Lin.\ Alg.\ Appl.} {\bf 23} (1979), 141--156.

\bibitem{Ray}
Ray, D.\ \ 
Sojourn times of diffusion processes.  
\emph{Illinois J.\ Math.} {\bf 7} (1963), 615--630. 

\end{thebibliography}
\end{document}